\documentclass[11pt,twoside,reqno,centertags,draft]{amsart}
\usepackage{amsfonts}
\usepackage{color,enumitem,graphicx}
\usepackage[colorlinks=true,urlcolor=blue,
citecolor=red,linkcolor=blue,linktocpage,pdfpagelabels,
bookmarksnumbered,bookmarksopen]{hyperref}
  \usepackage{amsmath,amsthm,amsfonts,amssymb}

\begin{document}
\title{ The Dirichlet elliptic  problem  involving regional fractional Laplacian }
\date{}
\maketitle

\vspace{ -1\baselineskip}

{\small
\begin{center}

\medskip

 {\sc Huyuan Chen \footnote{chenhuyuan@yeah.net} %\quad{and}\quad  Ying Wang\footnote{yingwang00@126.com}
  }
 \bigskip

{\small
Department of Mathematics, Jiangxi Normal University,\\
Nanchang, Jiangxi 330022, PR China
}\\[1mm]

\end{center}
}

\renewcommand{\thefootnote}{}
%\footnote{Accepted for publication: January 2004.}
%\footnote{*Corresponding author.}
\footnote{AMS Subject Classifications: 35R11, 35D30, 35J61. }
\footnote{Key words: Regional Fractional Laplacian,  Green Kernel, Weak solution. }

\begin{quote}
%\footnotesize
{\bf Abstract.} In this paper, we study the solutions for
elliptic equations involving  regional fractional Laplacian
\begin{equation}\label{0}
  \arraycolsep=1pt\left\{
\begin{array}{lll}
 \displaystyle  (-\Delta)^\alpha_\Omega   u=f \qquad & {\rm in}\quad   \Omega,\\[1mm]
\phantom{ (-\Delta)^\alpha     }
 \displaystyle   u=g\quad & {\rm on}\quad   \partial  \Omega,
\end{array}\right.
\end{equation}
where $\Omega$ is a bounded open domain in $\mathbb{R}^N$ ($N\ge 2$) with $C^2$ boundary $\partial\Omega$,
 $\alpha\in(\frac12,1)$ and the operator $(-\Delta)^\alpha_\Omega$ denotes the regional fractional Laplacian.
We prove that when $g\equiv0$, problem (\ref{0}) admits a unique weak solution in the cases that $f\in L^2(\Omega)$, $f\in L^1(\Omega, \rho^\beta
dx)$ and $f\in \mathcal{M}(\Omega,\rho^\beta)$, here $\rho(x)={\rm dist}(x,\partial\Omega)$, $\beta=2\alpha-1$ and $\mathcal{M}(\Omega,\rho^\beta)$ is a space of all Radon measures $\nu$ satisfying $\int_\Omega \rho^\beta d|\nu|<+\infty.$
Finally, we  provide an Integral by Parts Formula for the classical solution of (\ref{0}) with general boundary data $g$.

\end{quote}

\newcommand{\N}{\mathbb{N}}
\newcommand{\R}{\mathbb{R}}
\newcommand{\Z}{\mathbb{Z}}

\newcommand{\cA}{{\mathcal A}}
\newcommand{\cB}{{\mathcal B}}
\newcommand{\cC}{{\mathcal C}}
\newcommand{\cD}{{\mathcal D}}
\newcommand{\cE}{{\mathcal E}}
\newcommand{\cF}{{\mathcal F}}
\newcommand{\cG}{{\mathcal G}}
\newcommand{\cH}{{\mathcal H}}
\newcommand{\cI}{{\mathcal I}}
\newcommand{\cJ}{{\mathcal J}}
\newcommand{\cK}{{\mathcal K}}
\newcommand{\cL}{{\mathcal L}}
\newcommand{\cM}{{\mathcal M}}
\newcommand{\cN}{{\mathcal N}}
\newcommand{\cO}{{\mathcal O}}
\newcommand{\cP}{{\mathcal P}}
\newcommand{\cQ}{{\mathcal Q}}
\newcommand{\cR}{{\mathcal R}}
\newcommand{\cS}{{\mathcal S}}
\newcommand{\cT}{{\mathcal T}}
\newcommand{\cU}{{\mathcal U}}
\newcommand{\cV}{{\mathcal V}}
\newcommand{\cW}{{\mathcal W}}
\newcommand{\cX}{{\mathcal X}}
\newcommand{\cY}{{\mathcal Y}}
\newcommand{\cZ}{{\mathcal Z}}

\newcommand{\abs}[1]{\lvert#1\rvert}
\newcommand{\xabs}[1]{\left\lvert#1\right\rvert}
\newcommand{\norm}[1]{\lVert#1\rVert}

\newcommand{\loc}{\mathrm{loc}}
\newcommand{\p}{\partial}
\newcommand{\h}{\hskip 5mm}
\newcommand{\ti}{\widetilde}
\newcommand{\D}{\Delta}
\newcommand{\e}{\epsilon}
\newcommand{\bs}{\backslash}
\newcommand{\ep}{\emptyset}
\newcommand{\su}{\subset}
\newcommand{\ds}{\displaystyle}
\newcommand{\ld}{\lambda}
\newcommand{\vp}{\varphi}
\newcommand{\wpp}{W_0^{1,\ p}(\Omega)}
\newcommand{\ino}{\int_\Omega}
\newcommand{\bo}{\overline{\Omega}}
\newcommand{\ccc}{\cC_0^1(\bo)}
\newcommand{\iii}{\opint_{D_1}D_i}

\numberwithin{equation}{section}

\vskip 0.2cm \arraycolsep1.5pt
\newtheorem{lemma}{Lemma}[section]
\newtheorem{theorem}{Theorem}[section]
\newtheorem{definition}{Definition}[section]
\newtheorem{proposition}{Proposition}[section]
\newtheorem{remark}{Remark}[section]
\newtheorem{corollary}{Corollary}[section]

\setcounter{equation}{0}
\section{Introduction}

 The usual Laplaciain operator may be thought as a macroscopic manifestation of the Brownian motion, as known from
the Fokker-Plank equation for a stochastic
differential equation with a  Brownian motion (a Gaussian process), whereas the fractional Laplacian
operator $(-\Delta)^{ \alpha }$  is associated with
  a $2\alpha$-stable L\'{e}vy motion
(a non-Gaussian process) $L_t^{2\alpha}$, $\alpha\in(0,1)$, (see \cite{Duan2015} for a discussion about this microscopic-macroscopic relation).
From the observations and experiments related to L\'evy flights (\cite{BBC,KSZ,MK,SWS}), the fractional Laplacian described that a particle could  have infinite jumps
in an arbitrary time interval with intensity proportional to $\frac1{|x-y|^{N+2\alpha}}$, but if   the  particle jumping is forced to restrict only from one point $x\in\Omega$,
a  bounded open domain $\Omega$  in $\R^N$, to another point $y\in\Omega$ with the same intensity, then the related process is called the censored stable process and its generator is
the regional fractional Laplacian defined in $\Omega$,  see the references \cite{CK,CKS,GM}. In particular, the authors in \cite{BBC} pointed out that
the censored $2\alpha-$stable process is conservative
and will never approach $\partial\Omega$ when $\alpha\in (0,\frac12]$ and for $\alpha\in(\frac12,1)$ that process could approach to the boundary $\partial\Omega$. This indicates that
the Dirichlet problem involving the regional fractional Laplacian is well defined for $\alpha\in(\frac12,1)$ and in this note, we pay our attentions on  the solutions to related Dirichlet elliptic problem  with $\alpha\in(\frac12,1)$.

Throughout this paper, we assume that   $\alpha\in(\frac12,1)$, $\beta=2\alpha-1$, $\Omega$ is a bounded open domain in $\R^N$ ($N\ge 2$) with $C^2$ boundary $\partial\Omega$
and $\rho(x)={\rm dist}(x,\partial\Omega)$.
Denote by $(-\Delta)^\alpha_\Omega$   the regional fractional Laplacian
$$ (-\Delta)^\alpha_\Omega  u(x)=\lim_{\varepsilon\to0^+} (-\Delta)_{\Omega,\varepsilon}^\alpha u(x)$$
with
$$
(-\Delta)_{\Omega,\varepsilon}^\alpha  u(x)=-c_{N,\alpha}\int_{\Omega\setminus B_\varepsilon(x)}\frac{ u(z)-
u(x)}{|z-x|^{N+2\alpha}}\,dz,
$$
where $c_{N,\alpha}>0$ coincides the normalized constant of the fractional Laplacian.
The main objective of this note is to study the weak solution of   elliptic problem
\begin{equation}\label{eq 1.1}
  \arraycolsep=1pt\left\{
\begin{array}{lll}
 \displaystyle  (-\Delta)^\alpha_\Omega   u=f\qquad & {\rm in}\quad   \Omega,\\[1mm]
\phantom{ (-\Delta)^\alpha  }
 \displaystyle   u=0\quad & {\rm on}\quad   \partial  \Omega,
\end{array}\right.
\end{equation}
where $f:\Omega\to\R$. We will concentrate on the existence and uniqueness of the solution of (\ref{eq 1.1}) in a suitable weak sense when $f\in L^2(\Omega)$ or  $f$ belongs to Radon measure space.

When $f\in L^2(\Omega)$, it involves the Hilbert space $H_0^\alpha(\Omega)$ with
 the scalar product
$$\langle u,v\rangle_{\alpha}=\frac{c_{N,\alpha}}2\int_\Omega\int_\Omega \frac{[u(x)-u(y)][v(x)-v(y)] }{|x-y|^{N+2\alpha}}\,dxdy+\int_\Omega uv \,dx,\quad \forall\, u,v\in H_0^\alpha(\Omega).$$
which is the closure of $C^2_c(\bar\Omega)$ under the norm
$$\norm{u}_{H^\alpha(\Omega)}=\sqrt{\langle u,u\rangle_{ \alpha}},$$
which, shown in \cite{EGE}, is equivalent to the Gagliardo norm  $\norm{u}_{H^\alpha(\Omega)}$ in  $H_0^\alpha(\Omega)$
$$\norm{u}_{H^\alpha_0(\Omega)}:=\left(\frac{c_{N,\alpha}}2\int_\Omega\int_\Omega \frac{[u(x)-u(y)]^2 }{|x-y|^{N+2\alpha}}\,dxdy\right)^{\frac12}$$
 and its scalar product of $\norm{\cdot}_{H_0^\alpha(\Omega)}$ is
$$\langle u,v\rangle_{H_0^\alpha(\Omega)}=\frac{c_{N,\alpha}}2\int_\Omega\int_\Omega \frac{[u(x)-u(y)][v(x)-v(y)] }{|x-y|^{N+2\alpha}}\,dxdy,\quad \forall\, u,v\in H_0^\alpha(\Omega).$$

\begin{definition}
$(i)$  When $f\in L^2(\Omega)$,  a function  $u\in H^\alpha_0(\Omega)$,  is a weak solution of (\ref{eq 1.1}), if
$$
\langle u,v\rangle_{H_0^\alpha(\Omega)}=\int_\Omega f(x)v(x)\,dx,\qquad \forall\, v\in H^\alpha_0(\Omega).
$$

$(ii)$ %When $f\in L^1(\Omega,\rho^\beta \,dx)$, a function  $u\in L^1(\Omega)$ is a very weak solution of (\ref{eq 1.1}), if
%$$  \int_\Omega u(-\Delta)^\alpha_\Omega\xi \,dx=\int_\Omega f \xi \,dx,\quad \forall\xi\in \mathbb{X}_\alpha(\Omega),$$
%where $\rho(x)={\rm dist}(x,\partial\Omega)$ and
Denote by $\mathbb{X}_{\alpha}$ the space of functions
$\xi$,  continuous up to the boundary,  taking zero value on $\partial \Omega$ and verifying
$$\norm{(-\Delta)^\alpha_\Omega\xi}_{L^\infty(\Omega)}<\infty,$$
%$(iii)$
and by $\mathcal{M}(\Omega,\rho^\beta)$  the space of all the Radon measure $\nu$ satisfying
$$\int_\Omega \rho^\beta\, d|\nu|<+\infty.$$
When $f\in \mathcal{M}(\Omega,\rho^\beta)$, a function  $u\in L^1(\Omega)$ is a very weak solution of (\ref{eq 1.1}), if
$$
 \int_\Omega u(-\Delta)^\alpha_\Omega\xi \,dx=\int_\Omega \xi \,df,\quad \forall\xi\in \mathbb{X}_\alpha(\Omega).
$$

\end{definition}

We notice that $\beta=1$ if $\alpha=1$, and in this case the test functions' space corresponding to very weak solution is  $C^{1,1}_0(\Omega)$,
in which the function could be controlled by the distance function $\rho$. In the regional fractional case, the test functions space
$\mathbb{X}_\alpha(\Omega)$ plays the same role and the function in $\mathbb{X}_\alpha(\Omega)$ has the decay $\rho^\beta$, see Lemma \ref{lm 2.0} below.

Now we are ready to state our main theorem on the existence   and uniqueness of weak solution for problem (\ref{eq 1.1}).

\begin{theorem}\label{teo 1}
%Assume that  $\alpha\in(\frac12,1)$, $\beta=2\alpha-1$, $\Omega$ is a bounded open domain in $\R^N$ ($N\ge 2$) with $C^2$ boundary $\partial\Omega$.

$(i)$ Let $f\in L^2(\Omega)$,  then problem (\ref{eq 1.1}) has a unique weak  solution $u_f$ such that
\begin{equation}\label{3.002}
 \norm{u_f}_{H^\alpha_0(\Omega)}\le c_1\norm{f}_{L^2(\Omega)},
\end{equation}
where $c_1>0$.
\smallskip

$(ii)$ %Let $f\in L^1(\Omega,\rho^\beta \,dx)$, then problem (\ref{eq 1.1}) has a unique very weak  solution $u_f$ such that
%\begin{equation}\label{1.1}
%\norm{u_f}_{L^1(\Omega)}\le c_1\norm{f}_{L^1(\Omega,\rho^\beta \,dx)},
%\end{equation}
%where $c_1>0$. %$(iii)$
Let  $f\in \mathcal{M}(\Omega,\rho^\beta)$, then problem (\ref{eq 1.1}) has a unique very weak  solution $u_f$ such that
\begin{equation}\label{1.5}
\norm{u_f}_{L^1(\Omega)}\le c_2\norm{f}_{\mathcal{M}(\Omega,\rho^\beta)},
\end{equation}
where $c_2>0$.

\end{theorem}

   For $f\in L^2(\Omega)$ or $f\in \mathcal{M}(\Omega,\rho^\beta)$,    a sequence of functions  $\{f_n\}_n$ in $C^2(\Omega)\cap C(\bar \Omega)$ could be chosen to converge to $f$ in $L^2(\Omega)$ and we prove the solution of (\ref{eq 1.1}) is approached by the classical solution of
$$\arraycolsep=1pt\left\{
\begin{array}{lll}
 \displaystyle  (-\Delta)^\alpha_\Omega   u=f_n\qquad & {\rm in}\quad   \Omega,\\[2mm]
\phantom{ (-\Delta)^\alpha  }
 \displaystyle   u=0\quad & {\rm on}\quad   \partial  \Omega.
\end{array}\right.
$$
In this approaching process,   the most important  tool is  the Integral by Parts formula,
\begin{equation}\label{1.4}
\arraycolsep=1pt
\begin{array}{lll}
  \displaystyle\int_\Omega u(-\Delta)^\alpha_\Omega v \,dx&=& \displaystyle\frac{c_{N,\alpha}}2\int_\Omega\int_\Omega \frac{[u(x)-u(y)][v(x)-v(y)] }{|x-y|^{N+2\alpha}}\,dxdy
 \\[4mm]&=& \displaystyle\int_\Omega v(-\Delta)^\alpha_\Omega u \,dx, \qquad \forall u,v\in \mathbb{X}_\alpha(\Omega).
\end{array}
\end{equation}
 Thanks to a fractional Hardy-Sobolev inequality from \cite{D}, we also show the equivalence between  the norms $\norm{u}_{H^\alpha(\Omega)}$ and $\norm{\cdot}_{H_0^\alpha(\Omega)}$ for functions in $C_0^\infty(\Omega)$.

   It is known that  $L^1(\Omega,\rho^\beta \,dx)$ is a proper subset of $\mathcal{M}(\Omega,\rho^\beta)$ and we abuse the notation without confusion that 
   $df(x)=f(x) dx$ when $f\in L^1(\Omega,\rho^\beta \,dx)$ in the definition of very weak solution.
 But the proofs of the existence
of very weak solutions to (\ref{eq 1.1}) are very different for $f $ in $L^1(\Omega,\rho^\beta \,dx)$ and in $\mathcal{M}(\Omega,\rho^\beta)$. For $f\in L^1(\Omega,\rho^\beta \,dx)$, the very weak solution is approached directly by a Cauchy sequence in $L^1(\Omega)$, while in the case of $f\in \mathcal{M}(\Omega,\rho^\beta)$, we have to prove the approximations
is uniformly bounded in $L^1(\Omega)$ and uniformly integrable, then Dunford-Pettis Thoerem is applied to derive the very weak solution of (\ref{eq 1.1}). The elliptic problems involving
measure data  with second order operators have been extensively studied in \cite{BB1,BP,GV,P,V} and the reference therein,  and recently,
the elliptic problems involving the fractional Laplacian have been investigated by \cite{CV1,CV2,CV3}.

Finally, we make use of the nonlocal characteristic property to build an Integral by Parts Formula for the solution $u$ of
\begin{equation}\label{eq 1.2}
  \arraycolsep=1pt\left\{
\begin{array}{lll}
 \displaystyle  (-\Delta)^\alpha_\Omega   u=f\quad & {\rm in}\quad   \Omega,\\[1mm]
\phantom{ (-\Delta)^\alpha  }
 \displaystyle   u=g\quad & {\rm on}\quad   \partial  \Omega
\end{array}\right.
\end{equation}
for $f\in C^2( \Omega)\cap C(\bar \Omega)$ and $g\in C^2(\partial\Omega)$. That is,
\begin{equation}\label{1.3}
 \int_\Omega u(-\Delta)^\alpha_\Omega\xi \,dx=\int_\Omega f(x)\xi(x)\,dx+\int_{\partial \Omega}\frac{\partial^{\beta}\xi(x)}{\partial \vec{n}_x^{\beta}}g(x)\,d\omega(x),\quad \forall\xi\in \mathbb{X}_\alpha(\Omega)\cap \mathbb{D}_\beta,
 \end{equation}
where  $\vec{n}_x$ is the unit exterior normal vector  of $\Omega$ at point $x\in\partial\Omega$ and
\begin{equation}\label{4.2}
 \mathbb{D}_\beta=\bigcup_{ \tau\ge \beta}\{\phi_1 \rho^\tau+\phi_2:\ \phi_1,\phi_2\in C^2(\bar\Omega) \}.
\end{equation}
Here we note that
$$\frac{\partial^{\beta}\xi(x)}{\partial \vec{n}_x^{\beta}}=\lim_{t\to0^+}\frac{\xi(x)-\xi(x-t\vec{n}_x)}{t^\beta}=-\lim_{t\to0^+}\frac{\xi(x-t\vec{n}_x)}{t^\beta}$$
and an Integral by Parts Formula provided in \cite{G} states as follows
\begin{equation}\label{1.2}
   \int_\Omega w(-\Delta)^\alpha_\Omega v \,dx  =   \int_\Omega v (-\Delta)^\alpha_\Omega w  \,dx, \quad \forall w,v\in \mathbb{D}_\beta.
 \end{equation}

The paper is organized as follows. In Section \S2, we  study the solutions of (\ref{eq 1.1}) with $f\in L^\infty$, including  the existence
and uniqueness of classical solution, the boundary regularities and also provides some important estimates for proving (\ref{1.4}).   Section  \S3 is devoted
to give an Integration by Parts Formula for $u,v\in \mathbb{X}_\alpha(\Omega)$, then we obtain the existence and uniqueness of weak solution of (\ref{eq 1.1}) with zero boundary data
and $f\in L^2(\Omega)$. In Section \S4, we prove the very weak solution of (\ref{eq 1.1}) for $f\in L^1(\Omega,\rho^\beta \,dx)$ and $f\in \mathcal{M}(\Omega,\rho^\beta)$.
Finally, we provide  Integral by Parts Formula (\ref{1.3}) for the solution of (\ref{eq 1.1}) with general boundary data in Section \S5.

\setcounter{equation}{0}
\section{Preliminary}
The purpose of this section is to introduce some preliminaries on the classical solution of (\ref{eq 1.1}). We start it by
 the Comparison Principle. In what follows,  we denote by $c_i$ a generic  positive constant.

\begin{lemma}\label{comparison}
Assume that $g$ is continuous on $\partial \Omega$ and  $f_i:\Omega\to\R$ with $i=1,2$ are
continuous functions satisfying
$$f_1\ge f_2\quad {\rm in}\quad \Omega.$$
Let $u_1$ and $u_2$ be the solutions  of (\ref{eq 1.2}) with $f=f_1$ and $f_2$, respectively.
 Then
\begin{equation}\label{2.01}
 u_1  \ge u_2 \quad{\rm in}\quad  \Omega.
\end{equation}
Furthermore, if $f\equiv0,\ g\equiv 0$, then problem (\ref{eq 1.2}) only has zero solution.
\end{lemma}
{\bf Proof.}  By contradiction, if (\ref{2.01}) fails, denoting  $w=u_1-u_2$,   there exists $x_0\in\Omega$ such that
$$w(x_0)=u_1(x_0)-u_2(x_0)=\min_{x\in\Omega}w(x)<0.$$
Combining with  $w=0$ on $\partial\Omega$, we observe that
$$(-\Delta)^{\alpha}_\Omega   w(x_0)=   \int_{\Omega } \frac{w(x_0)-w(y)}{|x_0-y|^{N+2\alpha}}\,dy<0.$$
But
$$ (-\Delta)^{\alpha}_\Omega w(x_0)=f_1(x_0)-f_2(x_0)\ge 0,$$
which is impossible.  \qquad $\Box$
\medskip

Our main aim  here is give the regularity up to the boundary of the solution of (\ref{eq 1.1}).

\begin{theorem}\label{teo 6.1}
Assume that  $\alpha\in(\frac12,1)$,   $f\in L^\infty(\Omega)$,  $f_\pm=\max\{\pm f,0\}$ and $\rho(x)=dist(x,\partial\Omega)$.
Then problem (\ref{eq 1.1})  admits a unique solution $u_f$ such that
\begin{equation}\label{6.1.3}
-c_3 \norm{f_-}_{L^\infty(\Omega)}  \rho(x)^{\beta}\le  u_f(x)\le c_3\norm{f_+}_{L^\infty(\Omega)}\rho(x)^{\beta}\quad x\in\Omega.
\end{equation}
Moreover,

$(i)$ for $\theta\in (0,2\alpha)$ and an open set $\mathcal{O}$ satisfying $d_{\mathcal{O}}:={\rm dist}(\mathcal{O},\partial \Omega)>0$,
 there exists $c_4>0$ dependent of $d_{\mathcal{O}}$ and $\theta$ such that
\begin{equation}\label{6.1.2}
 \norm{u_f}_{C^\theta(\mathcal{O})}\le c_4 \norm{f}_{L^{\infty}(\Omega)};
\end{equation}

$(ii)$   there exists $c_5>0$
independent of $f$ such that
\begin{equation}\label{6.1.1}
\norm{u_f}_{C^{\beta}(\bar\Omega)}\le c_5\norm{f}_{L^\infty(\Omega)}.
\end{equation}

\end{theorem}

In order to consider (\ref{eq 1.1}), we need the
following uniformly estimates. Denote by $G_{\Omega,\alpha}$ the Green kernel of $(-\Delta)^\alpha_\Omega$ in $\Omega\times\Omega$ and by $\mathbb{G}_{\Omega,\alpha}[\cdot]$ the
Green operator defined as
$$\mathbb{G}_{\Omega,\alpha}[f](x)=\int_{\Omega} G_{\Omega,\alpha}(x,y)f(y)\,dy. $$

\begin{lemma}\label{lm 6.2.1}
Let  $\alpha\in(\frac12,1)$ and  $f\in L^\infty(\Omega)$,
 then $\mathbb{G}_{\Omega,\alpha}[f]$ is the unique solution of problem (\ref{eq 1.1})   and
\begin{equation}\label{annex 01q1}
\left|\mathbb{G}_{\Omega,\alpha}[f](x)\right|\le c_6\rho(x)^{\beta},\quad\forall\, x\in\Omega.
\end{equation}

\end{lemma}
 {\bf Proof.} The uniqueness follows by Lemma \ref{comparison}.  We observe that $\mathbb{G}_{\Omega,\alpha}[f]$ is a solution of  $ (-\Delta)^\alpha_\Omega w=f$ in $\Omega$
 and   for any $(x,y)\in
\Omega\times\Omega$ with $x\neq y$,
\begin{equation}\label{annex 6.01}
G_{\Omega,\alpha}(x,y)\le c_7
\min\left\{\frac1{|x-y|^{N-2\alpha}},\, \frac{\rho(x)^{\beta}\rho(y)^{\beta}}{|x-y|^{N-1+\beta}} \right\},
\end{equation}
see \cite{CK}. Then we have that
\begin{eqnarray*}
 | \mathbb{G}_{\Omega,\alpha}[f](x)| &\le & c_7\int_\Omega \frac{\rho(x)^{\beta}\rho(y)^{\beta}}{|x-y|^{N-1+\beta}}  |f(y)|dy\\
    &\le & c_7\rho(x)^{\beta}\norm{f}_{L^\infty (\Omega)}\int_\Omega  \frac{\rho(y)^{\beta}}{|x-y|^{N-1+\beta}}dy\\
   &\le& c_8\norm{f}_{L^\infty (\Omega)}\rho(x)^{\beta}, \quad \forall\,x\in\Omega.
\end{eqnarray*}
Hence $\mathbb{G}_{\Omega,\alpha}[f]$ is a solution of (\ref{eq 1.1})   verifying  (\ref{annex 01q1}).\qquad$\Box$\smallskip

In what follows, we  denote
 $$u_f=\mathbb{G}_{\Omega,\alpha}[f].$$

\begin{lemma}\label{lm 6.3.1}

For any $x_0\in \Omega$ and $\theta\in(0,2\alpha)$, there exists $c_9>0$ independent of $\rho(x_0)$ such that
\begin{equation}\label{6.3.1}
\norm{u_f}_{C^\theta(B_{\rho_0}(x_0))}\le c_9\rho_0^{\beta-\theta}  \norm{f}_{L^\infty(\Omega)},
\end{equation}
where $\rho_0=\rho(x_0)/3$.

\end{lemma}
 {\bf Proof.} For $x_0\in\Omega$, we denote that $\Omega_0=\{y\in \R^N:\ x_0+\rho_0y\in\Omega \}$ and
 $$v_f(x)= u_f(x_0+\rho_0 x),\quad \forall x\in \R^N,$$
 then by Lemma \ref{lm 6.2.1}, we have that
 $$\norm{v_f}_{L^\infty(B_2(0))}=\norm{u_f}_{L^\infty(B_{2\rho_0}(x_0))}\le c_6\norm{f}_{L^\infty(\Omega)}\rho_0^{\beta} $$
and for $x\in B_2(0)$,
\begin{eqnarray*}
 (-\Delta)^\alpha_{\Omega_0} v_f(x) &=&{\rm P.V.}c_{N,\alpha}\int_{\Omega_0} \frac{ u_f(x_0+\rho_0 x)- u_f(x_0+\rho_0 y)}{|x-y|^{N+2\alpha}}dy   \\
    &=&\rho_0^{2\alpha} (-\Delta)^\alpha_{\Omega} u_f(x_0+\rho_0x)
    = \rho_0^{2\alpha}f(x_0+\rho_0x)
\end{eqnarray*}
and
\begin{eqnarray*}
 (-\Delta)^\alpha  v_f(x)  &=&  (-\Delta)^\alpha_{\Omega_0}  v_f(x)+ v_f(x)\phi_0(x)\\
 &=& \rho_0^{2\alpha}f(x_0+\rho_0x) + v_f(x)\phi_0(x),\quad \forall x\in B_2(0),
\end{eqnarray*}
where $\phi_0(x)=c_{N,\alpha}\int_{\R^N\setminus \Omega_0} \frac1{|x-y|^{N+2\alpha}}dy.$
Since $B_3(0)\subset \Omega_0$, we have that
$$\phi_0(x)\le   c_{N,\alpha},\quad \forall x\in B_2(0).$$
Then by \cite[Proposition 2.3]{RS},
we have that
 \begin{eqnarray*}
 \norm{v_f}_{C^\gamma(B_1(0))} &\le& c_{10} \left( \norm{\rho_0^{2\alpha}f(x_0+\rho_0\cdot)+ v_f\phi_0}_{L^\infty(B_2(0))}+\norm{v_f}_{L^\infty(B_2(0))}\right)\\
   &\le &   c_{10} \left(\rho_0^{2\alpha}\norm{f}_{L^\infty(\Omega)}+\norm{v_f}_{L^\infty(B_2(0))}\right).
 \end{eqnarray*}
Since
$$\norm{v_f}_{L^\infty(B_2(0))}=\norm{u_f}_{L^\infty(B_{2\rho_0}(x_0))}\le c_{11}\rho_0^{\beta}\norm{f}_{L^\infty(\Omega)},$$
then we have that
$$\norm{u_f}_{ C^\theta(B_{\rho_0}(x_0))}\le c_{12}  \rho_0^{\beta-\theta}\norm{f}_{L^\infty(\Omega)}. $$
The proof ends.\qquad$\Box$

\medskip
 \begin{lemma}\label{lm 6.2.2}
 Let
\begin{equation}\label{con 1}
 \phi(x)=\int_{\R^N\setminus \Omega} \frac1{|x-y|^{N+2\alpha}}dy,
\end{equation}  then  $\phi\in C^{0,1}_{\rm loc}(\Omega)$ and for some $c_{13}>1$
\begin{equation}\label{6.2.6}
 \frac1{c_{13}}\rho(x)^{-2\alpha}\le \phi(x)\le c_{13}\rho(x)^{-2\alpha},\quad\forall x\in\Omega.
\end{equation}
\end{lemma}
{\bf Proof.} For
 $x_1,x_2\in \Omega $ and  any $z\in \R^N\setminus \Omega$, we have that
 $$|z-x_1|\ge \rho(x_1)+\rho(z), \qquad |z-x_2|\ge \rho(x_2)+\rho(z)$$ and
$$\left||z-x_1|^{N+2\alpha}-|z-x_2|^{N+2\alpha}\right|\leq
c_{14}|x_1-x_2|(|z-x_1|^{N+2\alpha-1}+|z-x_2|^{N+2\alpha-1}),
$$
for some $c_{14}>0$ independent of $x_1$ and $x_2$. Then
\begin{eqnarray*}
&&|\phi(x_1)-\phi(x_2)| \leq\int_{\R^N\setminus \Omega}  \frac{||z-x_2|^{N+2\alpha}-|z-x_1|^{N+2\alpha}|}{|z-x_1|^{N+2\alpha}|z-x_2|^{N+2\alpha}}dz
\\&&\leq c_{14}|x_1-x_2|\left[\int_{\R^N\setminus \Omega}  \frac{dz}{|z-x_1||z-x_2|^{N+2\alpha}}+\int_{\R^N\setminus \Omega}  \frac{dz}{|z-x_1|^{N+2\alpha}|z-x_2|}\right].
\end{eqnarray*}
By direct computation, we have that
\begin{eqnarray*}
 \int_{\R^N\setminus \Omega}  \frac{1}{|z-x_1||z-x_2|^{N+2\alpha}}dz
  &\le& \int_{\R^N\setminus{B_{\rho(x_1)}(x_1)}} \frac{1}{|z-x_1|^{N+2\alpha+1}}dz
 \\&&+\int_{\R^N\setminus{B_{\rho(x_2)}(x_2)}} \frac{1}{|z-x_2|^{N+2\alpha+1}}dz \\
 &\le& c_{15}[\rho(x_1)^{-1-2\alpha}+\rho(x_2)^{-1-2\alpha}]
\end{eqnarray*}
and similar to obtain that
$$\int_{\R^N\setminus \Omega} \frac{1}{|z-x_1|^{N+2\alpha}|z-x_2|}dz\le  c_{15}[\rho(x_1)^{-1-2\alpha}+\rho(x_2)^{-1-2\alpha}].$$
%where $c_9>0$ is independent of $x_1,x_2$.
Then
$$|\phi(x_1)-\phi(x_2)|\le c_{14} c_{15}[\rho(x_1)^{-1-2\alpha}+\rho(x_2)^{-1-2\alpha}]|x_1-x_2|,$$
that is,
 $\phi$ is $C^{0,1}$ locally in $\Omega$.

We next prove (\ref{6.2.6}). Without loss of generality, we may assume that
$0\in\partial\Omega$, the inside pointing normal vector at $0$ is $e_N=(0,\cdots,0,1)\in\R^N$ and %$x=se_N$ with
let $s\in(0,\frac14)$ such that $\R^N\setminus \Omega\subset \R^N\setminus B_s(se_N)$ and for $t>0$, we denote the cone
$$A_t=\{y=(y',y_N)\in\R^N: y_N\le s-t|y'|\}.$$
We observe that
there is $c_0>0$ such that
$$\left[A_{t_0}\cap \left(B_1(se_N)\setminus B_{2s}(se_N)\right)\right]\subset \R^N\setminus \Omega. $$
By the definition of $\phi$, we have that
 \begin{eqnarray*}
\phi(se_N)= \int_{\R^N\setminus \Omega} \frac1{|se_N-y|^{N+2\alpha}}dy \le  \int_{\R^N\setminus B_s(se_N)} \frac1{|se_N-y|^{N+2\alpha}}dy  \le c_{16}s^{-2\alpha}.
 \end{eqnarray*}
 On the other hand, we have that
\begin{eqnarray*}
 \int_{\R^N\setminus \Omega} \frac1{|se_N-y|^{N+2\alpha}}dy \ge  \int_{A_{c_0}\cap \left(B_1(se_N)\setminus B_{2s}(se_N)\right)} \frac1{|se_N-y|^{N+2\alpha}}dy
     \ge c_{17}s^{-2\alpha}.
 \end{eqnarray*}
The proof ends.\qquad $\Box$ \medskip

Now we are ready to prove Theorem \ref{teo 6.1}.
\medskip

\noindent {\bf Proof of Theorem \ref{teo 6.1}.}
 By Lemma \ref{lm 6.2.1}, $u_f$ is the unique solution of  (\ref{eq 1.1}).
 Since $\mathbb{G}_{\Omega,\alpha}[f_+]$, $-\mathbb{G}_{\Omega,\alpha}[f_-]$ are   solutions of  (\ref{eq 1.1}) replaced $f$ by $f_+$ and $f_-$ respectively.
 Then (\ref{6.1.3}) follows by Lemma \ref{comparison}.\smallskip

{\it Proof of (\ref{6.1.2}).}
Let $\tilde w=w$ in $\Omega$, $\tilde w=0$ in $\R^N\setminus\bar \Omega$,  we observe that
\begin{eqnarray*}
 (-\Delta)^\alpha \tilde w(x)  =  (-\Delta)^\alpha_\Omega w(x)+ w(x)\phi(x),\quad \forall\, x\in\Omega,
\end{eqnarray*}
where $\phi$ is defined by (\ref{con 1}).
It follows by Lemma \ref{lm 6.2.2}, $\phi\in C^{0,1}_{\rm loc}(\Omega)$.
Then
\begin{eqnarray*}
 (-\Delta)^\alpha \tilde w(x)  =f(x)+  w(x)\phi(x),\quad \forall\, x\in\Omega.
\end{eqnarray*}
Let $\mathcal{O}_1$ be a $C^2$ open set such that
$$\mathcal{O}\subset \mathcal{O}_1,\quad {\rm dist}(\mathcal{O}_1,\partial\Omega)= d_{\mathcal{O}}/2,\quad
{\rm and\quad dist}(\mathcal{O},\partial \mathcal{O}_1)= d_{\mathcal{O}}/2.$$
By \cite[Lemma 3.1]{CV3}, for any $\theta\in(0,2\alpha)$, we have that
\begin{equation}\label{6.2.4}
\norm{\tilde w}_{C^\theta(\mathcal{O})} \le  c_{18}\left[\norm{\tilde w}_{L^\infty(\mathcal{O}_1)}+\norm{\tilde w}_{L^1(\Omega)}+\norm{f+\tilde w\phi}_{L^\infty(\mathcal{O}_1)} \right].
\end{equation}
By Lemma \ref{lm 6.2.1} and Lemma \ref{lm 6.2.2}, we obtain that
$$\norm{\tilde w}_{L^\infty(\mathcal{O}_1)}+\norm{\tilde w}_{L^1(\Omega)}\le \norm{f}_{L^\infty(\Omega)}$$
and
$$|\tilde w(x)|\phi(x)\le c_{19}\rho(x)^{2\alpha-1}\norm{f}_{L^\infty(\Omega)}\rho(x)^{-2\alpha}\le c_{19}\rho(x)^{-1}\norm{f}_{L^\infty(\Omega)},$$
then
\begin{eqnarray*}
 \norm{f+\tilde w\phi}_{L^\infty(\mathcal{O}_1)} \le  \norm{f}_{L^\infty(\Omega)}+ \norm{\tilde w\phi}_{L^\infty(\mathcal{O}_1)}
   \le   c_{20}d_{\mathcal{O}}^{-1}\norm{f}_{L^\infty(\Omega)}.
\end{eqnarray*}
 Then (\ref{6.1.2}) holds.\smallskip

{\it Proof of (\ref{6.1.1}).}
  Taking $\theta=2\alpha-1$ in Lemma \ref{lm 6.3.1}, we have that
\begin{equation}\label{eq53151}
\frac{u(x)-u(y)}{|x-y|^{\theta}}\le c_{21}\|f\|_{L^\infty(\Omega)}
\end{equation}
for all $x,y$ such that $y\in B_R(x)$ with $R=\rho(x)/3$. We next show that (\ref{eq53151}) holds for all $x,y\in\bar\Omega$ with some renewed constant.

Indeed, we observe that after a Lipschitz change of coordinates, the bound   (\ref{eq53151}) remains the same except for the value of the constant $c$.
Then we can flatten the boundary near $x_0\in\partial\Omega$ to assume that $\Omega\cap B_{\rho_0}(x_0)=\{x_n>0\}\cap B_1(0)$.
Thus,  (\ref{eq53151}) holds for all $x,y$ satisfying $|x-y|\le \gamma x_n$
for some $\gamma=\gamma(\Omega)\in(0,1)$ dependent of the Lipschitz mapping.

Let $z=(z',z_n)$ and $w=(w',w_n)$ be two points in $\{x_n>0\}\cap B_{1/4}(0)$
and $r=|z-w|$. Denote that $\bar z=(z',z_n+r)$, $\bar w=(w',w_n+r)$ and $z_k=(1-\gamma^k)z+\gamma^k\bar z$
and $w_k=\gamma^k w+(1-\gamma^k)\bar w$, $k\ge0$. Then, using the bound  (\ref{eq53151}) whenever $|x-y|\le\gamma x_n$, we have that
$$|u(z_{k+1})-u(z_k)|\le c_{22}|z_{k+1}-z_k|^\theta=c_{21}|\gamma^k(z-\bar z)(\gamma-1)|^\theta\le c_{22}\gamma^k|z-\bar z|.$$
Moreover, since $x_n>r$ in all the segment joining $\bar z$ and $\bar w$, splitting this segment into a bounded number of segments
of length less than $\gamma r$, we obtain that
$$|u(\bar z)-u(\bar w)|\le c_{23}|\bar z-\bar w|^\theta\le c_{23}r^\theta.$$
Therefore,
 \begin{eqnarray*}
 |u(z)-u(w)|&\le&\sum_{k\ge0}|u(z_{k+1})-u(z_k)|+|u(\bar z)-u(\bar w)|+\sum_{k\ge0}|u(w_{k+1})-u(w_k)|
 \\&\le & (c_{24}\sum_{k\ge0}(\gamma^kr)^\theta+c_{25}r^\theta)(\|u\|_{L^\infty(\R^N)}+\|g\|_{L^\infty(\Omega)})
 \\&\le& c_{26}(\|u\|_{L^\infty(\R^N)}+\|g\|_{L^\infty(\Omega)})|z-w|^\theta,
 \end{eqnarray*}
which ends the proof.  \qquad$\Box$\smallskip

For a unbounded nonhomogeneous term $f$, we have that
\begin{lemma}\label{lm 2.0}
Assume that  $f$ is a $C^\gamma_{loc}(\Omega)$ function satisfying
$$|f(x)|\le c_{27}\rho(x)^{-\beta},\quad \forall x\in\Omega,$$
where $\gamma\in(0,1)$. Then problem (\ref{eq 1.1}) has a unique solution $u_f$ satisfying
\begin{equation}\label{2.0}
 |u_f(x)|\le c_{28}\norm{f\rho^{\beta}}_{L^\infty (\Omega)}\rho^\beta(x),\quad \forall x\in\Omega.
\end{equation}
%where $c_{28}>0$  depending on $\norm{f\rho^{\beta}}_{L^\infty (\Omega)}$.

\end{lemma}
{\bf Proof.} The uniqueness follows by Lemma \ref{comparison}.  It is known that $\mathbb{G}_{\Omega,\alpha}[f]$ is a solution of  $ (-\Delta)^\alpha_\Omega w=f$ in $\Omega$.
From (\ref{annex 6.01}), we have that for $x\in\Omega$,
\begin{eqnarray*}
 | \mathbb{G}_{\Omega,\alpha}[f](x)| &\le & c_7\int_\Omega \frac{\rho(x)^\beta\rho(y)^\beta}{|x-y|^{N-2+2\alpha}}  |f(y)|\,dy\\
    &\le & c_7\rho(x)^{\beta}\norm{f\rho^{\beta}}_{L^\infty (\Omega)}\int_\Omega  \frac{1}{|x-y|^{N-2+2\alpha}}\,dy\\
   &\le& c_7\norm{f\rho^{\beta}}_{L^\infty (\Omega)}\rho(x)^{\beta}\int_{B_{d_0}(x)} \frac{1}{|x-y|^{N-2+2\alpha}}\,dy ,
\end{eqnarray*}
where   $d_0=\sup_{x,y\in\Omega}|x-y|$ and $\int_{B_{d_0}(x)} \frac{1}{|x-y|^{N-2+2\alpha}}\,dy<+\infty$ by the fact that $N-2+2\alpha<N$.
Therefore, we obtain that $\mathbb{G}_{\Omega,\alpha}[f]$ is a solution of (\ref{eq 1.1}) satisfying (\ref{2.0}).\qquad$\Box$

\begin{remark}\label{re 2.1}
We remark that (\ref{2.0}) holds for $v\in \mathbb{X}_\alpha(\Omega)$.
In fact, let $f=(-\Delta)^\alpha_\Omega v$, which
satisfies
$$\norm{f\rho^{\beta}}_{L^\infty (\Omega)}<+\infty.$$
%Then    (\ref{2.0}) holds  for $v\in \mathbb{X}_\alpha(\Omega)$.
\end{remark}

\medskip

The next proposition plays an important role in the proof of Integration by Parts Formula with nonzero Dirichlet boundary condition.
For this purpose, we introduce some notations.
Denote
\begin{equation}\label{2.10}
 \Omega_\delta:=\{x\in\Omega:\ \rho(x)>\delta\}\quad{\rm and}\quad A_\delta:=\{x\in\Omega:\ \rho(x)<\delta\}.
\end{equation}
Since $\Omega$ is $C^2$, there exists $\delta_0>0$ such that
$\Omega_{\delta}$ is $C^2$ for any $\delta\in(0,\delta_0]$ and it is known that
for any $x\in \partial \Omega_{\delta}$, there exists $x^*\in \partial\Omega$ such that
$$|x-x^*|=\rho(x)\quad {\rm and}\quad x=x^*+\rho(x)\vec{n}_{x^*}.$$

\begin{proposition}\label{pr 2.2}
Assume that   $f\in C^2( \Omega)\cap C(\bar \Omega)$
and $g\in C^2(\partial\Omega)$.
Let $u$ be the classical solution of (\ref{eq 1.2}).
Then $u\in C^2(\Omega) \cap C^\beta(\bar\Omega)$.
Furthermore,  for $\delta\in(0,\delta_0)$, there exists $c_{29}>0$ such that
\begin{equation}\label{2.2}
|u(x)-u(y)|<c_{29}\rho(x)^{\beta-1}|x-y|,\quad \forall\, x\in \Omega_\delta,\ \forall\,  y\in A_\delta.
\end{equation}
\end{proposition}
{\bf Proof.}  {\it To prove $u\in C^2(\Omega) \cap C^\beta(\bar\Omega)$.} Here we only have to prove $u\in C^2(\Omega) \cap C^\beta(\bar\Omega)$
in the case that $g\equiv0$. In fact, since $\Omega$ is $C^2$ and $g\in C^2(\partial\Omega)$, then there exists $G\in C^2(\bar\Omega)$ such that
$$G=g\quad{\rm on}\ \partial\Omega.$$
Now we only consider the regularity of $u-G$, which is the solution of
\begin{equation}\label{eq 3.2}
  \arraycolsep=1pt
\begin{array}{lll}
 \displaystyle  (-\Delta)^\alpha_\Omega   u=f-(-\Delta)^\alpha_\Omega G\qquad & {\rm in}\quad   \Omega,\\[2mm]
\phantom{ (-\Delta)^\alpha  }
 \displaystyle   u=0\quad & {\rm on}\quad   \partial  \Omega.
\end{array}
\end{equation}
So it follows by Theorem \ref{teo 6.1}  that $u\in C^\beta(\bar \Omega)$.

We next prove $u\in C^2(\Omega)$. Extend the function $u$ by zero in $\R^N\setminus\Omega$, still denote it $u$,  and then
$$
   (-\Delta)^\alpha u(x)= (-\Delta)^\alpha_\Omega u(x)+u(x)\phi(x)=f(x)+u(x)\phi(x),\quad\forall x\in\Omega,
$$
where $(-\Delta)^\alpha$ is the global fractional Laplacian.
From Lemma \ref{lm 6.2.2}, $\phi\in C^{0,1}_{loc}(\Omega)$, applying \cite[Corollary 1.6]{RS} with $\theta<1+2\alpha$,
we have that
$$u\in C^\theta_{loc}(\Omega).$$
This means $u\in C^2(\Omega)$ since we can choose $\theta>2$.

{\it To prove (\ref{2.2}).}
By the compactness of $\partial\Omega$, we only consider a point $x_0\in \partial\Omega$ and
for simplicity, we can assume that $x_0=0$. Let $x=t\vec{n}_0$ with $t\in(0,\delta)$ and $y\in B_{\frac\delta3}(\delta \vec{n}_0)\cap \Omega_\delta$, for any $t\in(0,\delta]$, there exists $n$ depending on $t$ such that
$$t\in \left(\frac \delta{3^{n+1}},\frac {\delta}{3^{n}}\right)$$
and then we choose $x_k=x+\frac{1}{3^{k}}(y-x)$, $k=0,1,\cdots, n$.  We observe that
$$\rho(x_k)\ge   \frac{\delta}{3^{k}},\quad k=0,1,\cdots, n.$$

From Lemma \ref{lm 6.3.1}, it follows that for $k=0,1,\cdots, n$
\begin{equation}\label{2.001}
\norm{w}_{C^{0,1}(  B_{ \rho(x_k)/2}(x_k))}\le c_{30}\rho(x_k)^{\beta-1}  \norm{f}_{L^\infty(\Omega)}\le c_{30}(\frac{\delta}{3^{k}})^{\beta-1}  \norm{f}_{L^\infty(\Omega)}.
\end{equation}

  It is obvious that $x_0=y$ and
$$|x_k-x_{k+1}|<\frac{|x-y|}{3^k}.$$
Then we have that for $k=0,1,\cdots,n$,
\begin{eqnarray*}
 |w(x_k)-w(x_{k-1})| &\le & \norm{w}_{C^{0,1}( B_{\frac{\rho(x_k)}2}(x_k))}|x_k-x_{k-1}|  \\
   &\le &  c_{30} (\frac{\delta}{3^k})^{\beta-1} |x_k-x_{k-1}|
   \\
   &\le &  c_{30} t^{\beta-1} |x_k-x_{k-1}|,
\end{eqnarray*}
therefore,
\begin{eqnarray*}
 |w(x)-w(x_{0})|  \le   c_7 t^ {(1-\beta)} \sum_{k=0}^n \frac1{3^k} |x-y|
     \le   c_{31} t^{\beta-1}  |x-y|,
\end{eqnarray*}
where $c_{31}>0$ is independent of $t$.
So for some $c_{32}>0$, we have that
 \begin{eqnarray}\label{2.010}
 |w(x)-w(y)| \le    c_{32} \rho(x)^{\beta-1}  |x-y|.
\end{eqnarray}
For $y\in \Omega_\delta\setminus B_{\frac\delta3}(\delta \vec{n}_0)$, we may choose $y'\in B_{\frac\delta3}(\delta \vec{n}_0)\cap \Omega_\delta$.
There are at most $N_0$ points $y_k\in \Omega_\delta$ connecting $y$ and $y'$ such that
$$\frac\delta3\le |y_k-y_{k-1}|\le \frac\delta2.$$
We see that
\begin{eqnarray*}
 |w(y)-w(y')| \le    c_{33} \delta^{\beta-1}  |y-y'|.
\end{eqnarray*}
From (\ref{2.010}), we see that
\begin{eqnarray*}
 |w(x)-w(y')| \le    c_{32} \rho(x)^{\beta-1}  |x-y'|.
\end{eqnarray*}
Since $|y-y'|\ge \frac\delta3$ and $|x-y|>\delta$, then
\begin{eqnarray*}
  |w(x)-w(y)|&\le & |w(x)-w(y')|+ |w(y)-w(y')|\\
   &\le &  c_{33} \delta^{\beta-1}  |y-y'|+ c_9 \rho(x)^{\beta-1}  |x-y'|
   \\   &\le &c_{34}\rho(x)^{\beta-1} |x-y|.
\end{eqnarray*}
We finish the proof.  \qquad$\Box$

\begin{lemma}\label{lm 2.1}
Assume that    $f\in C^2( \Omega)\cap C(\bar \Omega)$, $g\in C^2(\partial\Omega)$  and  $w$ is  the classical solution of
(\ref{eq 1.2}).
Then
\begin{equation}\label{2.1}
\int_\Omega\int_\Omega\frac{[u(x)-u(y)]^2}{|x-y|^{N+2\alpha}}\,dxdy<+\infty.
\end{equation}

\end{lemma}
\noindent{\bf Proof. } From the interior regularity, we know that $u\in C^2(\Omega)\cap C^\beta(\bar\Omega)$.
From   \cite[Theorem 3.4]{GM}, it infers that
\begin{eqnarray}
 &&\frac{c_{N,\alpha}}2 \int_{\Omega_\delta}\int_{\Omega_\delta}\frac{[u(x)-u(y)]^2}{|x-y|^{N+2\alpha}}\,dxdy  \nonumber
 \\&&\qquad = \int_{\Omega_\delta}u(x)(-\Delta)^\alpha_{\Omega_\delta} u(x) \,dx \nonumber
 \\&&\qquad =  c_{N,\alpha} \int_{\Omega_\delta}\int_{ A_\delta}\frac{u(x)-u(y)}{|x-y|^{N+2\alpha}} u(x)dy \,dx+\int_{\Omega_\delta}u(x)(-\Delta)^\alpha_{\Omega} u(x) \,dx\nonumber
\\&&\qquad =  c_{N,\alpha} \int_{\Omega_\delta}\int_{ A_\delta}\frac{u(x)-u(y)}{|x-y|^{N+2\alpha}} u(x)dy \,dx+\int_{\Omega_\delta}u(x) f(x) \,dx. \label{2.3}
\end{eqnarray}
We observe that
\begin{equation}\label{2.4}
 |\int_{\Omega_\delta}u(x) f(x) \,dx|\le |\Omega|\norm{u}_{L^\infty(\Omega)} \norm{f}_{L^\infty(\Omega)}.
\end{equation}
From Proposition \ref{pr 2.2}, we derive that
\begin{eqnarray*}
   \int_{\Omega_\delta}\int_{ A_\delta}\frac{|u(x)-u(y)|}{|x-y|^{N+2\alpha}} |u(x)|dy \,dx&\le & \norm{u}_{L^\infty(\Omega)} \int_{ A_\delta} \int_{\Omega_\delta}\frac{|u(x)-u(y)|}{|x-y|^{N+2\alpha}} \,dx  dy \nonumber
 \\& \le& c_{35}\norm{u}_{L^\infty(\Omega)} \int_{ A_\delta}\rho(y)^{\beta-1} \int_{\Omega_\delta}\frac{1}{|x-y|^{N+2\alpha-1}} \,dx  dy\nonumber
 \\& \le& c_{36}\norm{u}_{L^\infty(\Omega)}  \int_{ A_\delta}\rho(y)^{\beta-1} \int_{\delta-\rho(y)}^{d_0}\frac{1}{r^{ 2\alpha}} dr  dy\nonumber
  \\& \le&c_{37}\norm{u}_{L^\infty(\Omega)}  \int_{ A_\delta}\rho(y)^{\beta-1} (\delta-\rho(y))^{-\beta} dy.
\end{eqnarray*}
Since $\Omega$ is $C^2$, then for $t\in(0,\delta)$ and $\delta\le \delta_0$, we have that
$$\frac12|\partial \Omega|\le |\partial \Omega_t|\le 2|\partial \Omega|$$
and by Fubini's theorem
\begin{eqnarray*}
\int_{ A_\delta}\rho(y)^{\beta-1} (\delta-\rho(y))^{-\beta} dy &=&  \int_0^{\delta} t^{\beta-1} (\delta-t)^{-\beta} |\partial \Omega_t| dt \\
    &\le & 2|\partial \Omega| \int_0^{\delta} t^{\beta-1} (\delta-t)^{-\beta}dt
   \\ &= & 2|\partial \Omega|\int_0^1 t^{\beta-1} (1-t)^{-\beta}dt.
\end{eqnarray*}
Therefore, for some $c_{38}>0$ independent of $\delta$ there holds that
 \begin{equation}\label{2.6}
\int_{\Omega_\delta}\int_{ A_\delta}\frac{|u(x)-u(y)|}{|x-y|^{N+2\alpha}}| u(x)|\,dy dx<c_{38},
 \end{equation}
thus, together with (\ref{2.3})-(\ref{2.4}), we derive that
$$\int_\Omega\int_ \Omega\frac{[u(x)-u(y)]^2}{|x-y|^{N+2\alpha}}\,dxdy<+\infty. $$
The proof ends. \qquad$\Box$

\begin{corollary}\label{cr 2.1}
Assume that    $f,h \in C^2( \Omega)\cap C(\bar \Omega)$ and  $u,w$ are the classical solution of
(\ref{eq 1.1}) with nonhomogeneous nonlinearities $f$ and $h$, respectively.

Then
\begin{equation}\label{2.5}
\lim_{\delta\to0^+} \int_{\Omega_\delta}\int_{ A_\delta}\frac{u(x)-u(y)}{|x-y|^{N+2\alpha}} w(y)dy \,dx =0.
\end{equation}

\end{corollary}
{\bf Proof.} From Theorem \ref{teo 6.1},
$$|w(x)|\le c_{15}\rho(x)^\beta,\qquad \forall\, x\in\Omega.$$
Thus,
\begin{eqnarray*}
   \int_{\Omega_\delta}\int_{ A_\delta}\frac{|u(x)-u(y)|}{|x-y|^{N+2\alpha}} |w(y)|dy \,dx&\le & \norm{w}_{L^\infty(A_\delta)} \int_{ A_\delta} \int_{\Omega_\delta}\frac{|u(x)-u(y)|}{|x-y|^{N+2\alpha}} \,dx  dy \nonumber
  \\& \le&c_{15}\delta^\beta \int_{ A_\delta}\rho^{\beta-1}(y) (\delta-\rho(y))^{-\beta} dy.
\end{eqnarray*}
By (\ref{2.6}), we have that
 $$\int_{\Omega_\delta}\int_{ A_\delta}\frac{|u(x)-u(y)|}{|x-y|^{N+2\alpha}} |w(y)|dy \,dx\le c_{38}\delta^\beta, $$
then (\ref{2.5}) holds. \qquad$\Box$

\setcounter{equation}{0}
\section{Zero boundary data}

\subsection{Classical solution}

In this subsection, we concentrate on the classical solution of (\ref{eq 1.1}) when $f\in C^2( \Omega)\cap C(\bar \Omega)$.

\begin{proposition}\label{pr 3.1}
Assume that   $f\in C^2( \Omega)\cap C(\bar \Omega)$  and  $u$ is the classical solution of
(\ref{eq 1.1}).
Then  
  \begin{equation}\label{3.1}
   \int_\Omega u(-\Delta)^\alpha_\Omega v \,dx  = \int_\Omega  f(x) v(x) \,dx,\quad \forall\, v\in \mathbb{X}_\alpha(\Omega).
 \end{equation}

\end{proposition}
{\bf Proof.} %From Lemma \ref{lm 2.1}, we deduce that $u$ satisfies (\ref{2.1}). We only have to prove (\ref{3.1}).
 Let $h(x)=(-\Delta)^\alpha_\Omega v(x)$ and $h_n$ be a sequence of  $C^2( \Omega)\cap C(\bar \Omega)$
functions such that
$$\lim_{n\to\infty}\norm{h_n-h}_{L^\infty(\Omega)}=0.$$
Let $v_n$ be the solution of
\begin{equation}\label{eq 3.1}
\arraycolsep=1pt\left\{
\begin{array}{lll}
 \displaystyle  (-\Delta)^\alpha_\Omega  u=h_n\qquad & {\rm in}\quad   \Omega,\\[2mm]
\phantom{ (-\Delta)^\alpha  }
 \displaystyle   u=0\quad & {\rm on}\quad   \partial  \Omega.
\end{array}\right.
\end{equation}
and then
$$\norm{v_n}_{C^{\beta}(\Omega)}\le c_{16}\norm{h}_{L^\infty(\Omega)}.$$
From Lemma \ref{lm 2.1}, we have that
$$\int_\Omega\int_\Omega\frac{[u(x)-u(y)]^2}{|x-y|^{N+2\alpha}}\,dxdy<+\infty\quad{\rm and}\quad \int_\Omega\int_\Omega\frac{[v_n(x)-v_n(y)]^2}{|x-y|^{N+2\alpha}}\,dxdy<+\infty,$$
which imply that
$$\int_{\Omega_\delta}\int_{\Omega_\delta}\frac{|[u(x)-u(y)][v_n(x)-v_n(y)]|}{|x-y|^{N+2\alpha}}\,dxdy<+\infty.$$

From \cite[Theorem 3.4]{GM}, it infers that
\begin{eqnarray*}
&&\frac{c_{N,\alpha}}2 \int_\Omega\int_\Omega\frac{[u(x)-u(y)][v_n(x)-v_n(y)]}{|x-y|^{N+2\alpha}}\,dxdy
 \\&&\qquad =\frac{c_{N,\alpha}}2 \lim_{\delta\to0^+} \int_{\Omega_\delta}\int_{\Omega_\delta}\frac{[u(x)-u(y)][v_n(x)-v_n(y)]}{|x-y|^{N+2\alpha}}\,dxdy  \nonumber
 \\&&\qquad =\lim_{\delta\to0^+} \int_{\Omega_\delta}v_n(x) (-\Delta)^\alpha_{\Omega_\delta} u(x) \,dx \nonumber
 \\&&\qquad =\int_{\Omega }v_n(x)(-\Delta)^\alpha_{\Omega} u(x) \,dx+  c_{N,\alpha}\lim_{\delta\to0^+} \int_{\Omega_\delta}\int_{ A_\delta}\frac{u(x)-u(y)}{|x-y|^{N+2\alpha}} v_n(x)\,dy dx\nonumber
\\&&\qquad = \int_{\Omega }v_n(x) f(x) \,dx+  c_{N,\alpha}\lim_{\delta\to0^+}  \int_{\Omega_\delta}\int_{ A_\delta}\frac{[u(x)-u(y)][v_n(x)-v_n(y)]}{|x-y|^{N+2\alpha}} \, dydx
\\&&\qquad\quad\ +c_{N,\alpha}\lim_{\delta\to0^+}  \int_{\Omega_\delta}\int_{ A_\delta}\frac{u(x)-u(y)}{|x-y|^{N+2\alpha}}v_n(y) \,dydx
\end{eqnarray*}
and by Corollary \ref{cr 2.1}, we have that
$$c_{N,\alpha}\lim_{\delta\to0^+}  \int_{\Omega_\delta}\int_{ A_\delta}\frac{[u(x)-u(y)][v_n(x)-v_n(y)]}{|x-y|^{N+2\alpha}} \,dydx=0$$
and
$$c_{N,\alpha}\lim_{\delta\to0^+}  \int_{\Omega_\delta}\int_{ A_\delta}\frac{u(x)-u(y)}{|x-y|^{N+2\alpha}}v_n(y) \,dydx=0.$$
Therefore,
\begin{equation}\label{3.02}
 \frac{c_{N,\alpha}}2 \int_\Omega\int_\Omega\frac{[u(x)-u(y)][v_n(x)-v_n(y)]}{|x-y|^{N+2\alpha}}\,dxdy=\int_{\Omega }v_n (-\Delta)^\alpha_{\Omega} u  \,dx=\int_{\Omega }v_n  f \,dx.
\end{equation}
Since $u$ and $v_n$ have the same role the above procedures, then
$$\frac{c_{N,\alpha}}2 \int_\Omega\int_\Omega\frac{[u(x)-u(y)][v_n(x)-v_n(y)]}{|x-y|^{N+2\alpha}}\,dxdy=\int_{\Omega } u (-\Delta)^\alpha_{\Omega} v_n  \,dx=\int_{\Omega }h_n  u  \,dx.$$
Therefore, (\ref{3.1}) holds.\qquad $\Box$

\medskip

From the above observations, we are ready to prove the Integral by Parts Formula for the regional fractional Laplacian.

\begin{theorem}\label{teo 3.1}
Let $u,v\in\mathbb{X}_\alpha(\Omega)$, then
\begin{equation}\label{3.0}
   \int_\Omega\int_\Omega \frac{[u(x)-u(y)]^2}{|x-y|^{N+2\alpha}}\,dxdy  \le  c_{17} \norm{(-\Delta)^\alpha_\Omega u}^2_{L^\infty(\Omega)}
 \end{equation}
 and
  \begin{equation}\label{3.01}
   \int_\Omega u(-\Delta)^\alpha_\Omega v \,dx=\frac{c_{N,\alpha}}2 \int_\Omega\int_\Omega \frac{[u(x)-u(y)][v(x)-v(y)] }{|x-y|^{N+2\alpha}}\,dxdy  =  \int_\Omega v(-\Delta)^\alpha_\Omega u \,dx.
 \end{equation}

\end{theorem}
{\bf Proof.}  Let $f(x)=(-\Delta)^\alpha_\Omega u(x)$, $h(x)=(-\Delta)^\alpha_\Omega v(x)$, and choose $\{f_n\}_n$, $\{h_n\}_n$ two sequences of  $C^2( \Omega)\cap C(\bar \Omega)$
functions such that
\begin{eqnarray}\label{3.05}
\lim_{n\to\infty}\norm{f_n-f}_{L^\infty(\Omega)}=0\quad{\rm and}\quad\lim_{n\to\infty}\norm{h_n-h}_{L^\infty(\Omega)}=0.
\end{eqnarray}
Let $u_n$ and $v_n$ be the solution of (\ref{eq 3.1}) with nonhomogeneous terms $f_n$ and $h_n$ respectively.
Integrating (\ref{eq 3.1}) with nonhomogeneous terms $f_n$ by $u_n$ and $v_n$ over $\Omega$, from (\ref{3.02}), we have that
\begin{eqnarray}\label{3.03}
\frac{c_{N,\alpha}}2\int_\Omega\int_\Omega \frac{[u_n(x)-u_n(y)]^2}{|x-y|^{N+2\alpha}}\,dxdy =\int_\Omega u_n(-\Delta)^\alpha_\Omega u_n\,dx= \int_\Omega u_n  f_n \,dx
\end{eqnarray}
and
\begin{eqnarray}\label{3.04}
\frac{c_{N,\alpha}}2\int_\Omega\int_\Omega \frac{[u_n(x)-u_n(y)][v_n(x)-v_n(y)]}{|x-y|^{N+2\alpha}}\,dxdy =\int_\Omega v_n(-\Delta)^\alpha_\Omega u_n\,dx= \int_\Omega v_n   f_n  \,dx.
\end{eqnarray}

Since $$\norm{u_n}_{C^\beta(\bar\Omega)}\le c_{18}\norm{f_n}_{L^\infty(\Omega)}\le c_{39}\norm{(-\Delta)^\alpha_\Omega u }_{L^\infty(\Omega)},$$
it infers from (\ref{3.03}) that
$$\int_\Omega\int_\Omega \frac{[u_n(x)-u_n(y)]^2}{|x-y|^{N+2\alpha}}\,dxdy \le c_{40} \norm{(-\Delta)^\alpha_\Omega u }^2_{L^\infty(\Omega)}.
$$
%where $c_{18},c_{19},c_{20}>0$.
This implies that for any $\epsilon>0$ and any $n\in\N$,
$$
\int_\Omega\int_\Omega \frac{[u_n(x)-u_n(y)]^2}{|x-y|^{N+2\alpha}}\chi_{(\epsilon,\infty)}(|x-y|)\,dxdy \le c_{41} \norm{(-\Delta)^\alpha_\Omega u }^2_{L^\infty(\Omega)},
$$
passing to the limit as $n\to\infty$, then we obtain that for any $\epsilon>0$,
$$
\int_\Omega\int_\Omega \frac{[u(x)-u(y)]^2}{|x-y|^{N+2\alpha}}\chi_{(\epsilon,\infty)}(|x-y|)\,dxdy \le c_{41} \norm{(-\Delta)^\alpha_\Omega u }^2_{L^\infty(\Omega)}.
$$
Since the left hand side of above inequality is decreasing with respective to $\epsilon>0$ and the right hand side is independent of $\epsilon$, so passing to the limit as  $\epsilon\to0^+$, we derive (\ref{3.0}).\smallskip

{\it To prove (\ref{3.01}).} It is obvious that $v$ verifies (\ref{3.0}).  Then
$\frac{u_n(x)-u_n(y)}{|x-y|^{\frac{N+2\alpha}2}} $ converges to $\frac{u(x)-u(y)}{|x-y|^{\frac{N+2\alpha}2}}$ in $L^2(\Omega\times\Omega)$
and $\frac{v_n(x)-v_n(y)}{|x-y|^{\frac{N+2\alpha}2}} $ converges to $\frac{v(x)-v(y)}{|x-y|^{\frac{N+2\alpha}2}}$ in $L^2(\Omega\times\Omega)$,
thus,
$$\lim_{n\to\infty}\int_\Omega\int_\Omega \frac{[u_n(x)-u_n(y)][v_n(x)-v_n(y)]}{|x-y|^{N+2\alpha}}\,dxdy=\int_\Omega\int_\Omega \frac{[u(x)-u(y)][v(x)-v(y)]}{|x-y|^{N+2\alpha}}\,dxdy,$$
which, together with (\ref{3.05}), implies that
$$\int_\Omega v(-\Delta)^\alpha_\Omega u \,dx= \frac{c_{N,\alpha}}2\int_\Omega\int_\Omega \frac{[u(x)-u(y)][v(x)-v(y)] }{|x-y|^{N+2\alpha}}\,dxdy.$$
The same to conclude that
$$\int_\Omega u(-\Delta)^\alpha_\Omega v \,dx= \frac{c_{N,\alpha}}2\int_\Omega\int_\Omega \frac{[u(x)-u(y)][v(x)-v(y)] }{|x-y|^{N+2\alpha}}\,dxdy$$
and (\ref{3.01}) holds.\qquad$\Box$

\subsection{Weak solution when $f\in L^2(\Omega)$ }

Our aim in this subsection is to consider the weak solution of (\ref{eq 1.1}) when the nonhomogeneous term $f$ is weakened from $L^\infty(\Omega)$ to $L^2(\Omega)$. To this end, we have to involve  the fractional Hilbert space
$H^\alpha_0(\Omega)$, which is the closure of $C^2_c(\Omega)$ under the norm of
\begin{equation}\label{3.2.1}
 \norm{u}_{H^\alpha(\Omega)}:=\left(\frac{c_{N,\alpha}}2\int_\Omega\int_\Omega \frac{[u(x)-u(y)]^2 }{|x-y|^{N+2\alpha}}\,dxdy\right)^{\frac12}+\norm{u}_{L^2(\Omega)}.
\end{equation}
This is  called as Gagliardo norm and we denote by $\norm{u}_{H^\alpha_0(\Omega)}$ the first part of (\ref{3.2.1}) on the right hand side,
which, we shall prove, is  a equivalent norm of $\norm{u}_{H^\alpha(\Omega)}$ in $H^\alpha_0(\Omega)$. Then we may say that the space $H_0^\alpha(\Omega)$ is the closure of $C^2_c(\bar\Omega)$ under the norm $\norm{\cdot}_{H^\alpha_0(\Omega)}$.
%Denote the scalar product
%$$\langle u,v\rangle_{H_0^\alpha(\Omega)}=\int_\Omega\int_\Omega \frac{[u(x)-u(y)][v(x)-v(y)] }{|x-y|^{N+2\alpha}}\,dxdy,\quad \forall u,v\in H_0^\alpha(\Omega).$$

We make use of  a Poincar\'{e} type inequality  to prove the equivalence of the norms $\norm{.}_{H^\alpha(\Omega)}$
and $\norm{\cdot}_{H^\alpha_0(\Omega)}$.
\begin{proposition}\label{pr 3.2}
 The norms
$\norm{\cdot}_{H^\alpha(\Omega)}$ and $\norm{\cdot}_{H^\alpha_0(\Omega)}$ are equivalent in $H_0^\alpha(\Omega)$.
\end{proposition}
{\bf Proof.} For $C^2$ bounded domain and $\alpha\in(\frac12,1)$,  it follows by \cite[Theorem 1.1]{D} that
$$\int_\Omega\frac{|u(x)|^2}{\rho^{2\alpha}(x)}\,dx\le c_{41}\int_\Omega\int_\Omega \frac{[u(x)-u(y)]^2 }{|x-y|^{N+2\alpha}}\,dxdy,\quad \forall u \in C_c^2(\Omega),$$
which implies that
\begin{equation}\label{3.2.2}
\norm{u}_{L^2(\Omega)}\le c_{42} \norm{u}_{H^\alpha_0(\Omega)},\qquad \forall u\in C_c^2(\Omega).
\end{equation}
Since $C^2_c(\Omega)$ is dense in $H^\alpha_0(\Omega)$, then   (\ref{3.2.2}) holds in $  H^\alpha_0(\Omega)$.
We omit the  left proof.    \qquad $\Box$\medskip

\noindent{\bf Proof of Theorem \ref{teo 1} part $(i)$.}  {\it Uniqueness.} Let $u,w$ be two  weak solutions of (\ref{eq 1.1}), then we derive that
$$
  \langle u-w, \, v \rangle_{H^\alpha_0(\Omega)} =0,\qquad \forall v\in  H^\alpha_0(\Omega).
$$
Taking $v=u-w\in H_0^\alpha(\Omega)$, we have that
 $$\norm{u-w}_{H^\alpha_0(\Omega)} =0.$$
Then we obtain the uniqueness.\smallskip

{\it Existence.} Let $\{f_n\}_n$ be a sequence of functions in $\mathbb{X}_\alpha(\Omega)$ satisfying
$$\lim_{n\to\infty}\norm{f_n-f}_{L^2(\Omega)}=0.$$
Let $u_n$ be the classical solution of (\ref{eq 1.1}) with nonhomogeneous term $f_n$.
Then
\begin{equation}\label{3.2.3}
\langle u_n,v\rangle_{H^\alpha_0(\Omega)}=\int_\Omega v f_n \,dx,\qquad \forall v\in H^\alpha_0(\Omega).
\end{equation}
From Theorem \ref{teo 3.1}, Proposition \ref{pr 3.2} and H\"{o}lder inequality,
we have that
\begin{eqnarray*}
\norm{u_n}^2_{H^\alpha_0(\Omega)} = \int_\Omega f_n u_n \,dx
    &\le & \norm{u_n} _{L^2(\Omega)}\norm{f_n} _{L^2(\Omega)}\\
    &\le & c_{43}\norm{u_n}_{H^\alpha_0(\Omega)} \norm{f_n} _{L^2(\Omega)}.
\end{eqnarray*}
 Then we have that
\begin{equation}\label{3.2.4}
 \norm{u_n}_{H^\alpha_0(\Omega)} \le c_{43}\norm{f_n}_{L^2(\Omega)}\le c_{44}\norm{f}^2_{L^2(\Omega)}.
\end{equation}
From \cite[Theorem 6.7, Theorem 7.1]{EGE}, the embedding: $H^\alpha_0(\Omega)\hookrightarrow L^2(\Omega)$ is  compact, then up to subsequence,  there exists $u\in L^2(\Omega)$ such that
$$u_n \to u\ \ {\rm in}\ \left(H^\alpha_0(\Omega)\right)'\quad {as}\ n\to \infty$$
and
$$\lim_{n\to\infty}\norm{u_n-u}_{L^2(\Omega)}=0.$$
Then from (\ref{3.2.3}),
we have that
$$\langle u,v\rangle_{H^\alpha_0(\Omega)}=\int_\Omega v f \,dx,\quad \forall\, v\in H^\alpha_0(\Omega),$$
that is,
$u$ is a weak solution of (\ref{eq 1.1}).  Taking $v=u$ above, we deduce (\ref{3.002}).\qquad $\Box$\bigskip

%{\it  To prove (\ref{3.002}). }
%From (\ref{3.2.4}), we have
%In fact, it infers from (\ref{3.2.4}) that for any $\epsilon>0$,
%$$
%\int_\Omega\int_\Omega \frac{[u_n(x)-u_n(y)]^2}{|x-y|^{N+2\alpha}}\chi_{(\epsilon,\infty)}(|x-y|)\,dxdy \le c \norm{f}^2_{L^2(\Omega)}
%$$
%and it could be written as
%$$
%\int_\Omega\int_\Omega \frac{u_n^2(x)+u_n^2(y)-2u_n(x)u_n(y)}{|x-y|^{N+2\alpha}}\chi_{(\epsilon,\infty)}(|x-y|)\,dxdy \le c \norm{f}^2_{L^2(\Omega)}
%$$ we have that
%$$\lim_{n\to\infty}\int_\Omega\int_\Omega \frac{u_n^2(x)}{|x-y|^{N+2\alpha}}\chi_{(\epsilon,\infty)}(|x-y|)\,dxdy
%=\int_\Omega\int_\Omega \frac{u^2(x)}{|x-y|^{N+2\alpha}}\chi_{(\epsilon,\infty)}(|x-y|)\,dxdy,$$
%$$\lim_{n\to\infty}\int_\Omega\int_\Omega \frac{u_n^2(y)}{|x-y|^{N+2\alpha}}\chi_{(\epsilon,\infty)}(|x-y|)\,dxdy
%=\int_\Omega\int_\Omega \frac{u^2(y)}{|x-y|^{N+2\alpha}}\chi_{(\epsilon,\infty)}(|x-y|)\,dxdy$$
%and
%$$\lim_{n\to\infty}\int_\Omega\int_\Omega \frac{u_n (x)u_n (y)}{|x-y|^{N+2\alpha}}\chi_{(\epsilon,\infty)}(|x-y|)\,dxdy
%=\int_\Omega\int_\Omega \frac{u (x)u(y)}{|x-y|^{N+2\alpha}}\chi_{(\epsilon,\infty)}(|x-y|)\,dxdy,$$
%so we deduce that
%$$
%\int_\Omega\int_\Omega \frac{[u(x)-u(y)]^2}{|x-y|^{N+2\alpha}}\chi_{(\epsilon,\infty)}(|x-y|)\,dxdy \le c \norm{f}^2_{L^2(\Omega)}
%$$
%which implies (\ref{3.002}).

\setcounter{equation}{0}
\section{Very weak solution }

\subsection{The case that $f\in L^1(\Omega,\rho^\beta \,dx)$ }
In this section, we may weaken the nonhomogeneous term $f$ to $L^1(\Omega,\rho^\beta \,dx)$ in (\ref{eq 1.1}).

\medskip

\noindent{\bf Proof of Theorem \ref{teo 1} part $(ii)$ when  $f\in L^1(\Omega,\rho^\beta \,dx)$.}  {\it Uniqueness.} Let $u,w$ be two very weak solutions of (\ref{eq 1.1}), then
\begin{equation}\label{3.3.0}
 \int_\Omega (u-w)(-\Delta)^\alpha_\Omega v \,dx =0,\qquad \forall v\in  \mathbb{X}_\alpha(\Omega).
\end{equation}
Let $\eta_{u-w}$ be the solution of
 \begin{equation}\label{test}
 \arraycolsep=1pt\left\{
\begin{array}{lll}
 \displaystyle  (-\Delta)^\alpha_\Omega   \eta_{u-w}={\rm sign}(u-w)\qquad & {\rm in}\quad   \Omega,\\[2mm]
\phantom{ (-\Delta)^\alpha  }
 \displaystyle  \eta_{u-w}=0\quad & {\rm on}\quad   \partial  \Omega.
\end{array}\right.
 \end{equation}
We observe that $\eta_{u-w}\in \mathbb{X}_\alpha(\Omega)$ and put $v=\eta_{u-w}$ in (\ref{3.3.0}), then we obtain that
$$\int_\Omega |u-w| \,dx =0,$$
which implies the uniqueness.\smallskip

{\it Existence.} We make use of Proposition \ref{pr 3.1} to approximate the solution $u$ of (\ref{eq 1.1}) by a sequence of classical solutions.
In fact, we choose a sequence of $ C^2( \Omega)\cap C(\bar \Omega)$
functions $\{f_n\}_n$ such that
\begin{equation}\label{3.3.1}
 \lim_{n\to\infty}\norm{f_n-f}_{L^1(\Omega,\,\rho^\beta \,dx)}=0.
\end{equation}
Denote $u_n$ the solution of (\ref{eq 1.1}) with nonhomogeneous nonlinearity $f_n$. Then
from Proposition \ref{pr 3.1}, we have that
 \begin{equation}\label{3.3}
   \int_\Omega u_n(-\Delta)^\alpha_\Omega v \,dx  = \int_\Omega  f_n(x) v(x) \,dx,\quad \forall v\in \mathbb{X}_\alpha(\Omega).
 \end{equation}
By Lemma \ref{lm 2.0}, it deduces that
\begin{equation}\label{4.1.1}
 |v(x)|\le c_{45}\rho^\beta(x),\quad \forall x\in \Omega
\end{equation}
and together with the convergence of $\{f_n\}_n$ in $L^1(\Omega,\rho^\beta \,dx)$, we obtain that
\begin{equation}\label{3.3.2}
 \lim_{n\to\infty}\int_\Omega f_n(x) v(x) \,dx=\int_\Omega f(x) v(x) \,dx.
\end{equation}

For any $n,m\in\N$, let $\eta_{u_m-u_n}$ be the solution of (\ref{test}) with nonhomogeneous term sign$(u_m-u_n)$,
then we obtain that
\begin{eqnarray*}
\int_\Omega |u_m-u_n| \,dx  =  \int_\Omega (f_m-f_n) \eta_{u_m-u_n} \,dx \le  c_{46}\int_\Omega |f_m-f_n|\rho^\beta \,dx.
\end{eqnarray*}
%where $c_{28}>0$.

For any $\epsilon>0$, it infers by (\ref{3.3.1}) that there exists $N_\epsilon>0$ such that for any $n,m>N_\epsilon$,
$$c_{47}\int_\Omega |f_m-f_n|\rho^\beta \,dx\le \epsilon,$$
which implies that for any $n,m>N_\epsilon$
  $$\int_\Omega |u_m-u_n| \,dx \le \epsilon.$$
Thus, $\{u_n\}_n$ is a Cauchy sequence in $L^1(\Omega)$ and then there exists $u\in L^1(\Omega)$ such that
$$\lim_{n\to\infty}\norm{u_n-u}_{L^1(\Omega)}=0.$$
Passing to the limit of (\ref{3.3}) as $n\to \infty$, we obtain that
$$ \int_\Omega u(x) (-\Delta)^\alpha_\Omega v(x) \,dx=\int_\Omega f(x) v(x) \,dx,\quad \forall\, v\in \mathbb{X}_\alpha(\Omega).$$
Therefore, problem (\ref{eq 1.1}) has a very weak solution, that is,
\begin{equation}\label{3.3.3}
  \int_\Omega u(-\Delta)^\alpha_\Omega v \,dx =\int_\Omega fv \,dx,\quad \forall\, v\in  \mathbb{X}_\alpha(\Omega),
\end{equation}
choosing $v=\eta_{u}$,  the solution of (\ref{test}) with nonhomogeneous term sign$(u)$, it infers from (\ref{3.3.3})
that
$$
  \int_\Omega |u|  \,dx \le c_{48}\int_\Omega |f|\rho^{\beta} \,dx.
$$
The proof ends. \qquad$\Box$

\subsection{The case that $f\in \mathcal{M}(\Omega,\rho^\beta)$ }
In this subsection, we may weaken the nonhomogeneous term $f$ to $\mathcal{M}(\Omega,\rho^\beta)$ in (\ref{eq 1.1}).\medskip

\noindent{\bf Proof of Theorem \ref{teo 1} part $(ii)$ when $f\in \mathcal{M}(\Omega,\rho^\beta)$.}  {\it Uniqueness.} Let $u,w$ be two very weak solutions of (\ref{eq 1.1}), then
$$
 \int_\Omega (u-w)(-\Delta)^\alpha_\Omega v \,dx =0,\quad \forall\, v\in  \mathbb{X}_\alpha(\Omega),
$$
which reduces to (\ref{3.3.0}).

{\it Existence.} We make use of Proposition \ref{pr 3.1} to approximate the solution $u$ of (\ref{eq 1.1}) by a sequence of classical solutions.
In fact, we choose a sequence of $ C^2( \Omega)\cap C(\bar \Omega)$
functions $\{f_n\}_n$ such that
\begin{equation}\label{4.3.1}
 \lim_{n\to\infty}\int_\Omega f_n v \,dx=\int_\Omega  v \,d f\quad{\rm for\ any}\quad v\in \mathbb{X}_\alpha(\Omega).
\end{equation}
Denote $u_n$ the solution of (\ref{eq 1.1}) with nonhomogeneous nonlinearity $f_n$. Then
from Proposition \ref{pr 3.1}, we have that
 \begin{equation}\label{4.3}
   \int_\Omega u_n(-\Delta)^\alpha_\Omega v \,dx  = \int_\Omega  f_n(x) v(x) \,dx,\quad \forall\, v\in \mathbb{X}_\alpha(\Omega).
 \end{equation}
Thus, it deduces by (\ref{4.1.1}) and  (\ref{4.3.1})  that
\begin{equation}\label{4.3.2}
 \lim_{n\to\infty}\int_\Omega v(x) f_n(x)\,dx=\int_\Omega  v(x) \,df(x).
\end{equation}

{\it To prove that $\{u_n\}_n$ is  uniformly bounded in $L^1(\Omega)$.} For any $n\in\N$, let $\eta_{u_n}$ be the solution of (\ref{test}) with nonhomogeneous term sign$(u_n)$,
then we derive that
\begin{eqnarray*}
\int_\Omega |u_n| \,dx  &=& \int_\Omega  f_n  \eta_{u_n} \,dx \le  c_{49}\int_\Omega  |\eta_{u_n}|\, d|f|
\\&\le& c_{49}c_{48}\norm{\eta_{u_n}\rho^{-\beta}}_{L^\infty(\Omega)}\int_\Omega  \rho^{\beta} d|f|\le c_{50}\int_\Omega  \rho^{\beta} d|f|,
\end{eqnarray*}
where $c_{49},c_{50}>0$ are independent of $n$, since $|\eta_{u_n}|\le c_{48}\rho^\beta$ in $\Omega$.

{\it To prove that $\{u_n\}_n$ is uniformly integrable.} Let $\mathcal{O}$ be any Borel subset of $\Omega$, take $\eta_{\mathcal{O}}$ be the solution of (\ref{test}) with nonhomogeneous term $\chi_{\mathcal{O}}{\rm sign}(\mathcal{O})$,
then we see that
\begin{eqnarray*}
\int_{\mathcal{O}} |u_n| \,dx  &=& \int_\Omega  f_n  \eta_{\mathcal{O}} \,dx \le  c_{49}\int_\Omega  |\eta_{\mathcal{O}}|\, d|f|
\\&\le& c_{49}c_{48}\norm{\eta_{\mathcal{O}}\rho^{-\beta}}_{L^\infty(\Omega)}\int_\Omega  \rho^{\beta} d|f|\le c_{51}\int_\Omega  \rho^{\beta} d|f|,
\end{eqnarray*}
where $c_{51}>0$ is independent of $n$.
We observe that
\begin{eqnarray*}
|\rho^{-\beta}(x)\eta_{\mathcal{O}}(x)| &= &\rho^{-\beta}(x)|\int_\Omega G_{\Omega,\alpha}(x,y) \chi_{\mathcal{O}}(y){\rm sign}(u_n)(y)dy|  \\
    &\le & c_7\int_{\mathcal{O}}\frac{\rho^\beta(y)}{|x-y|^{N-2+2\alpha}}dy
    \\
    &\le &c_7D_0^\beta\int_{B_{d_0}(x)}\frac1{|x-y|^{N-2+2\alpha}}dy
    \\ &\le & c_{52} d_0^{2-2\alpha}=  c_{52} |\mathcal{O}|^{\frac{2-2\alpha}{N}},
\end{eqnarray*}
where $c_{52}>0$ is independent of $n$, $D_0=\sup_{x,y\in\omega}|x-y|$ and $d_0>0$  satisfying
$$|\mathcal{O}|=|B_{d_0}(0)|.$$
 Thus, we derive that
\begin{eqnarray*}
\int_{\mathcal{O}} |u_n| \,dx  \le  c_{53}\norm{f}_{\mathcal{M}_\alpha(\Omega)}|\mathcal{O}|^{\frac{2-2\alpha}{N}},
\end{eqnarray*}
where $c_{53}>0$ is independent of $n$.

Therefore, we conclude that $\{u_n\}_n$ is uniformly bounded in $L^1(\Omega)$ and uniformly integrable, thus weakly compact in $L^1(\Omega)$ by the Dunford-Pettis Theorem,
and there exists a subsequence  $\{u_{n_k}\}_k$ and an integrable function $u$ such that $u_{n_k}\to u$ weakly in $L^1(\Omega)$. Passing to the limit in (\ref{4.3}), we obtain that
$u$ is a very weak solution of (\ref{eq 1.1}).

Choosing $v=\eta_{u}$,  the solution of (\ref{test}) with nonhomogeneous term sign$(u)$, it infers
that
$$
  \int_\Omega |u|  \,dx \le c_{54}\int_\Omega\rho^{\beta} d|f|.
$$
This ends the proof. \qquad$\Box$

\setcounter{equation}{0}
\section{General boundary data}

In this section, we consider the classical solution of (\ref{eq 1.2}) under the general boundary data.
In \cite[Theorem 1.3 $(i)$]{G}, the author proved the Integral by Part Formula
\begin{equation}\label{4.1}
   \int_\Omega u(-\Delta)^\alpha_\Omega v \,dx  =  \int_\Omega v(-\Delta)^\alpha_\Omega u \,dx+\int_{\partial\Omega}
   v\frac{\partial^\beta u}{\partial\vec{n}^\beta}\,d\omega-\int_{\partial\Omega}
   u\frac{\partial^\beta v}{\partial\vec{n}^\beta}\,d\omega,\quad \forall u,v\in  \mathbb{D}_\beta,
 \end{equation}
where $\mathbb{D}_\beta$ is given by (\ref{4.2}).

However, it is open to show that the solution $u_{f,g}$ of (\ref{eq 1.2})  belongs to
$ \mathbb{D}_\beta$, even under the hypothesis that
$$
  f\in C^2( \Omega)\cap C(\bar \Omega)\quad {\rm and}\quad g\in C^2(\partial\Omega).
$$

\begin{proposition}\label{pr 4.1}
Assume that   $f\in C^2( \Omega)\cap C(\bar \Omega)$, $g\in C^2(\partial\Omega)$   and  $u$ be the classical solution of
(\ref{eq 1.2}).
Then
  \begin{equation}\label{3.4}
   \int_\Omega u(-\Delta)^\alpha_\Omega v \,dx  = \int_\Omega  f  v  \,dx+\int_{\partial\Omega} g\frac{\partial^\beta v}{\partial \vec{n}^\beta} \,d\omega,\quad \forall v\in \mathbb{X}_\alpha(\Omega)\cap \mathbb{D}_\beta.
 \end{equation}
 \end{proposition}
\noindent {\bf Proof. } Since $\Omega$ is a $C^2$ domain and $g\in C^2(\partial\Omega)$, then there exists a $C^2(\bar\Omega)$ function $G$ such that
$$G=g\quad {\rm on}\quad \partial\Omega.$$

 Let $u_g$
 be the solution of (\ref{eq 1.2})  and denote
$$u_0=u_g-G.$$
Then
$u_0$ satisfies (\ref{eq 3.2}).
From \cite[Proposition 2.3]{GM}, $(-\Delta)^\alpha_\Omega G\in C^{2-2\alpha}_{loc}(\Omega)$ and
\begin{equation}\label{3.2}
|(-\Delta)^\alpha_\Omega G(x)|\le c_{55}\rho(x)^{-\beta},\qquad \forall x\in\Omega.
\end{equation}
Choose $\tilde g_n$ a sequence of $C^2( \Omega)\cap C(\bar \Omega)$ functions such that
$$ \lim_{n\to\infty}\norm{\left(\tilde g_n- (-\Delta)^\alpha_\Omega G\right)\rho^\beta}_{L^\infty(\Omega)}=0.$$
Let $w_n$ be the solution of
\begin{equation}\label{eq 3.4}
  \arraycolsep=1pt\left\{
\begin{array}{lll}
 \displaystyle  (-\Delta)^\alpha_\Omega   u=f-\tilde g_n\qquad & {\rm in}\quad   \Omega,\\[2mm]
\phantom{ (-\Delta)^\alpha  }
 \displaystyle   u=0\quad & {\rm on}\quad   \partial  \Omega.
\end{array}\right.
\end{equation}
By Proposition \ref{pr 3.1}, it infers that for $v\in \mathbb{X}_\alpha(\Omega)$,
\begin{eqnarray}
\int_\Omega w_n(x)(-\Delta)^\alpha_\Omega v(x)\,dx  &=&\int_\Omega v(x)(-\Delta)^\alpha_\Omega w_n(x)\,dx \nonumber    \\
    &=& \int_\Omega v(x)f(x)\,dx -\int_\Omega v(x) \tilde g_n\,dx.\label{3.5}
\end{eqnarray}
We observe that
$$\norm{w_n-u_g}_{L^\infty(\Omega)}\le c_{37}\norm{\left(\tilde g_n- (-\Delta)^\alpha_\Omega G\right)\rho^\beta}_{L^\infty(\Omega)}.$$
Therefore, passing to the limit of (\ref{3.5}) as $n\to\infty$, we have that
\begin{eqnarray*}
\int_\Omega u_0(x)(-\Delta)^\alpha_\Omega v(x)\,dx =  \int_\Omega v(x)f(x)\,dx  -\int_\Omega v(x)(-\Delta)^\alpha_\Omega G(x) \,dx.
\end{eqnarray*}
Since $G\in C^2(\bar \Omega)$, then for $v\in \mathbb{D}_\beta$, it deduces by (\ref{4.1}) that
$$\int_\Omega v(x)(-\Delta)^\alpha_\Omega G(x) \,dx=\int_\Omega G(x)(-\Delta)^\alpha_\Omega v(x)\,dx-\int_{\partial\Omega}g \frac{\partial^\beta v}{\partial \vec{n}^\beta}\,d\omega.$$
Thus,
$$\int_\Omega u_g(-\Delta)^\alpha_\Omega v \,dx=\int_\Omega (u_0+G)(-\Delta)^\alpha_\Omega v \,dx=\int_\Omega v f \,dx+\int_{\partial\Omega}g \frac{\partial^\beta v}{\partial \vec{n}^\beta}\,d\omega.$$
The proof ends.\qquad$\Box$\medskip

\begin{remark}
The function space $C_c^2(\Omega)$ is a subset of $\mathbb{X}_\alpha(\Omega)\cap \mathbb{D}_\beta$, so $\mathbb{X}_\alpha(\Omega)\cap \mathbb{D}_\beta$ is not empty.
However, $C_c^2(\Omega)$ is not proper for the Integral by Parts Formula since the boundary term
$\int_{\partial\Omega} g\frac{\partial^\beta v}{\partial \vec{n}^\beta} \,d\omega$ always vanishes for  $v\in C_c^2(\Omega)$.

\end{remark}

\bigskip\bigskip

\noindent{\bf Acknowledgements:}     H. Chen is supported by NSFC, No: 11401270,
and by the Jiangxi Provincial Natural Science Foundation, No: 20161ACB20007 and the
Project-sponsored by SRF for ROCS, SEM. H. Chen would like to thank professor Patricio Felmer for the supporting
in the visiting CMM of University of Chile.

\end{document}